\documentclass[11pt]{amsart}
\usepackage{amsfonts}
\usepackage{amssymb}
\usepackage[all]{xy}
\usepackage[british,french]{babel}
\usepackage[latin1]{inputenc}
\newtheorem{thm}{Theorem}[section]

\newtheorem{lemma}[thm]{Lemma}
\newtheorem{sublemma}[thm]{Sublemma}
\newtheorem{prop}[thm]{Proposition}
\theoremstyle{definition}
\newtheorem{defn}[thm]{Definition}
\theoremstyle{remark}
\newtheorem{rem}[thm]{Remark}

\usepackage{fancyhdr}
\numberwithin{equation}{section}

\newcommand{\gener}[1]{\langle #1 \rangle}
\newif\ifprivate
\privatefalse

\def\???{\ifprivate {\bf {???}} \marginpar{{\Huge {\bf ?}}}
\else \fi}

\newcommand{\RM}{\mathbb{R}}
\newcommand{\NM}{\mathbb{N}}
\newcommand{\C}{\mathcal{C}}
\newcommand{\F}{\mathcal{F}}
\newcommand{\Cinf}{\mathcal{C}^\infty}
\newcommand{\mybar}{\overline}
\newcommand{\phy}{\varphi}
\renewcommand{\O}{\mathcal{C}} % oui je sais c'est bizarre mais c'est
{\left\{\begin{array}{ll}} {\end{array}\right.}
\begin{document}
\title{ A singular Poincar\'{e} lemma} \author{Eva Miranda}
\address{Departament d'\`{A}lgebra i Geometria, Facultat de Matem\`{a}tiques, Universitat de Barcelona, Gran Via
 de les Corts Catalanes 585, 08007 Barcelona, Spain}
\email{evamiranda@ub.edu}
\thanks{ The first author is partially supported by the DGICYT project number BFM2003-03458. }

 \author{Vu Ngoc
  San}
  \address{Institut Fourier, Unit\'{e} mixte de recherche CNRS-UJF 5582, B.P. 74, 38402, Saint-Martin d'H\`{e}res, France}
\email{svungoc@ujf-grenoble.fr}
   \selectlanguage{british}
\begin{abstract}
  We prove a Poincar\'{e} lemma for a set of $r$ smooth functions on a
  $2n$-dimensional smooth manifold satisfying a commutation relation
  determined by $r$ singular vector fields associated to a Cartan subalgebra of $\mathfrak{sp}(2r,\mathbb
  R)$. This result has a natural interpretation in terms of the
  cohomology associated to the infinitesimal deformation of a
  completely integrable system.
\end{abstract}
\date{May 21, 2004}
% ----------------------------------------------------------------
\maketitle
\section{Introduction}
The classical Poincar\'{e} lemma asserts that a closed $1$-form on a
smooth manifold is locally exact. In other words, given $m$-functions
$g_i$ on an $m$-dimensional manifold for which
$\frac{\partial}{\partial x_i}(g_j)=\frac{\partial}{\partial
  x_j}(g_i)$ there exists a smooth $F$ in a neighbourhood of each
point such that $g_i=\frac{\partial}{\partial x_i}(F)$.

Now assume that we have a set of $r$ functions $g_i$ and a set of $r$
vector fields $X_i$ with a singularity at a point $p$ and fulfilling a
commutation relation of type $X_i(g_j)=X_j(g_i)$. We want to know if a
similar expression for $g_i$ exists in a neighbourhood of $p$.

In the case $g_i$ are $n$ functions on the symplectic manifold
$(\mathbb R^{2n},\sum_i dx_i\wedge dy_i)$ and $X_i$ form a basis of a
Cartan subalgebra of ${\frak{sp}}(2n,\mathbb R)$ a Poincar\'{e}-like
lemma exists. This result was stated by Eliasson in \cite{eli1}. In
\cite{eli2} Eliasson provided a proof of this statement in the
completely elliptic case. As far as the non-elliptic cases are
concerned, no complete proof of this result is known to the authors of
this note.

The analytical counterpart of this result dates back to the seventies
and was proved by Vey \cite{vey}. The transition from the analytical
case to the smooth case in cases other than elliptic entails a
non-trivial work with flat functions along certain submanifolds and,
in our opinion, cannot be neglected.

The aim of this note is to prove a more general singular Poincar\'{e}
lemma; the one that would correspond to a set of $r$ functions on a
$2n$-dimensional manifold with $r\leq n$ fulfilling similar
commutation relations determined by a basis of a Cartan subalgebra of
$\frak{sp}(2r,\mathbb R)$.  In particular, in this way we obtain a
complete proof also when $r=n$ in the non-completely elliptic cases
which was missing in the literature. This result has a natural
interpretation in terms of the cohomology associated to the
infinitesimal deformation of completely integrable foliations (see
section~\ref{sec:deformation}).

This result has applications in establishing normal forms for
completely integrable systems. The statement for $r=n$ was used by
Eliasson in \cite{eli1} and \cite{eli2} to give a symplectic normal
form for non-degenerate singularities of completely integrable
systems.  The more general result we prove here could be useful to
establish normal forms for more general singularities of completely
integrable systems.

\section{The result}

All the objects considered in this note will be $\Cinf$.  We are
interested in germ-like objects attached to a point $p$ of a smooth
manifold $M^{2n}$.

We denote by $(x_1,y_1,\dots,x_n,y_n)$ a set of coordinates centered
at the origin. Consider the standard symplectic form
$\omega=\sum_{i=1}^n dx_i\wedge dy_i$ in a neighbourhood of the
origin. Take $r\leq n$ and consider the embedding $i_{r}\colon\mathbb
R^{2r}\longrightarrow \mathbb R^{2n}$ defined by
$i_r(x_1,y_1,\dots,x_r,y_r)=(x_1,y_1,\dots,x_r,y_r,0,\dots,0)$.
Consider $\omega_r=\sum_{i=1}^r dx_i\wedge dy_i$ then
$i_r^*(\omega)=\omega_r$, in other words, this embedding induces an
inclusion of Lie groups $Sp(2r,\mathbb R)\subset Sp(2n,\mathbb R)$. In
this way $\frak{sp}(2r,\mathbb R)$ is realized as a subalgebra of
$\frak{sp}(2n,\mathbb R)$. This particular choice of subalgebra is
implicit throughout the note.

In this note we consider singular vector fields which constitute a
basis of a Cartan subalgebra of the Lie algebra $\frak{sp}(2r,\mathbb
R)$ with $r\leq n$. Recall that $\frak{sp}(2m,\mathbb R)$ is
isomorphic to the algebra of quadratic forms in $2m$ variables,
$Q(2m,\mathbb R)$, via symplectic duality. Thus the above chosen
immersion induces, in turn, an inclusion of subalgebras $Q(2r,\mathbb
R)\subset Q(2n,\mathbb R)$.

Cartan subalgebras of $Q(2r,\mathbb R)$ were classified by Williamson
in \cite{williamson}.

\begin{thm}{\label{Willi}}{\bf (Williamson)}
  For any Cartan subalgebra $\mathcal C$ of $Q(2r,\mathbb R)$ there is
  a symplectic system of coordinates $(x_1,y_1,\dots,x_r,y_r)$ in
  $\mathbb R^{2r}$ and a basis $q_1,\dots,q_r$ of $\mathcal C$ such
  that each $q_i$ is one of the following:
  \begin{equation}
    \begin{array}{lcr}
      q_i = x_i^2 + y_i^2 & \textup{for }  1 \leq i \leq k_e \ , &
\textup{(elliptic)}  \\
      q_i = x_iy_i  &  \textup{for }  k_e+1 \leq i \leq k_e+k_h \ , &
      \textup{(hyperbolic)}\\
      \begin{cases}
        q_i = x_i y_i + x_{i+1} y_{i+1} , \\
        q_{i+1} = x_i y_{i+1}-x_{i+1} y_i
      \end{cases} &
      \begin{array}{c}
        \textup{for }  i = k_e+k_h+ 2j-1,\\ 1\leq j \leq k_f
      \end{array}
      &
      \textup{(focus-focus pair)}
 \end{array}
  \end{equation}
\end{thm}

Observe that the number of elliptic components $k_e$, hyperbolic
components $k_h$ and focus-focus components $k_h$ is therefore an
invariant of the algebra $\mathcal C$. The triple $(k_e,k_h,k_f)$ is
called the Williamson type of $\mathcal C$. Observe that
$r=k_e+k_h+2k_f$.  Let $q_1,\dots,q_r$ be a Williamson basis of this
Cartan subalgebra. We denote by $X_i$ the Hamiltonian vector field of
$q_i$ with respect to $\omega$. Those vector fields are a basis of the
corresponding Cartan subalgebra of $\frak{sp}(2r,\mathbb R)$. We say
that a vector field $X_i$ is hyperbolic (resp. elliptic) if the
corresponding function $q_i$ is so. We say that a pair of vector
fields $X_i,X_{i+1}$ is a focus-focus pair if $X_i$ and $X_{i+1}$ are
the Hamiltonian vector fields associated to functions $q_i$ and
$q_{i+1}$ in a focus-focus pair.

In the local coordinates specified above, the vector fields $X_i$ take
the following form:

\begin{itemize}
\item $X_i$ is an elliptic vector field, $$X_i=
  2(-y_i\frac{\partial}{\partial x_i} +x_i{\frac{\partial}{\partial
      y_i}}).$$
\item $X_i$ is a hyperbolic vector field,
  $$X_i=-x_i\frac{\partial}{\partial x_i}+y_i{\frac{\partial}{\partial
      y_i}}.$$
\item $X_i,X_{i+1}$ is a focus-focus pair,
  $$X_i=-x_i\frac{\partial}{\partial
    x_{i}}+y_{i}\frac{\partial}{\partial
    y_i}-x_{i+1}\frac{\partial}{\partial
    x_{i+1}}+y_{i+1}\frac{\partial}{\partial y_{i+1}}$$
  and
  $$X_{i+1}=-x_i\frac{\partial}{\partial
    x_{i+1}}+y_{i+1}\frac{\partial}{\partial
    y_i}+x_{i+1}\frac{\partial}{\partial
    x_i}-y_i\frac{\partial}{\partial y_{i+1}}.$$
\end{itemize}

With all this notation at hand we can now state the result proven in
this note.

\begin{thm}
  \label{theo:principal}
  Let $g_1,\dots g_r$, be a set of germs of smooth functions on
  $(\RM^{2n},0)$ with $r\leq n$ fulfilling the following commutation
  relations
  $$X_i(g_j)=X_j(g_i), \quad \forall i,j \in\{1,\dots,r\}$$
  where the
  $X_i$'s are the vector fields defined above.  Then there exists a
  germ of smooth function $G$ and $r$ germs of smooth functions $f_i$
  such that,
\begin{enumerate}
\item $X_j(f_i)=0$, $\forall i,j \in\{1,\dots,r\}$.
\item $g_i=f_i+X_i(G)$ $\forall i \in\{1,\dots,r\}$.
\end{enumerate}
\end{thm}
\section{Preliminaries}
In this section we recall some basic facts which are proved elsewhere
and which will be used in the proof. Here and in the rest of the
article the symbol $X_i$ always refers to the hamiltonian vector field
associated to the quadratic function $q_i$, as precised above.

\subsection{A special decomposition for elliptic vector fields}

Assume $X_i$ is an elliptic vector field. That is, it is the vector
field associated to an elliptic $q_i=x_i^2+y_i^2$. The following
result was proved by Eliasson in \cite{eli1} when $n=1$.

\begin{prop}{\label{decomeli}}
  Let $g$ be a smooth function then there exist differentiable
  functions $g_1$ and $g_2$ such that:
  $$g=g_1(x_1,y_1,\dots,x_i^2+y_i^2,\dots, x_n,y_n)+X_i(g_2).$$
  Moreover,
  \begin{enumerate}
  \item $g_1$ is uniquely defined and satisfies $X_j(g_1)=0$ whenever
    $X_j(g)=0$;
  \item one can choose $g_2$ such that $X_j(g_2)=0$ whenever
    $X_j(g)=0$.
  \end{enumerate}
\end{prop}
\vspace{5mm}
\noindent\textbf{Remark: }
There are explicit formulas for the functions $g_1$ and $g_2$ claimed
above.  Let $\phi_t$ be the flow of the vector field $X_i$ we define,
$$g_1(x_1,y_1,\dots,x_n,y_n)=\frac{1}{\pi}\int_0^{\pi}
g(\phi_t(x_1,y_1,\dots,x_n,y_n))dt$$
and
$$g_2(x_1,y_1,\dots,x_n,y_n)=\frac{1}{\pi}\int_0^{\pi}
(tg(\phi_t(x_1,y_1\dots,x_n,y_n))-g_1(x_1,y_1\dots,x_n,y_n))dt.$$
\subsection{A special decomposition for hyperbolic vector fields}
In this section we assume the vector field $X_i$ corresponds to a
hyperbolic function $q_i=x_i y_i$. As a matter of notation, $S_i$
stands for the set $S_i=\{x_i=0,y_i=0\}$. When we refer to an
$(x_i,y_i)$-flat function f along $S_i$ we mean that
$$\frac{\partial^{k+l} f}{\partial x_i^k\partial
  y_i^l}_{\vert{S_i}}=0.$$
The first result is a decomposition result
for smooth functions.
\begin{prop}\label{hyp1}
  Given a smooth function $g$ there exist smooth functions $g_1$ and
  $g_2$ such that
  $$g=g_1(x_1,y_1,\dots, x_i.y_i, \dots,x_n,y_n)+X_i(g_2).$$
  Moreover
  one can choose $g_1$ and $g_2$ such that $X_j(g_1)=X_j(g_2)=0$
  whenever $X_j(g)=0$ for some $j\neq i$.
\end{prop}
\vspace{5mm} This proposition was proven by the first author of this
note in \cite{evathesis} (Proposition 2.2.2).

The main strategy of the proof is first to find a decomposition of
this type in terms of $(x_i,y_i)$-jets and then solve the similar
problem for $(x_i,y_i)$-flat functions along $S_i$. A main ingredient
in the proof of the proposition above are the following lemmas which
we will be also used in the proof of the theorem in this note. The
proof of the following two lemmas is also contained in
\cite{evathesis} (lemmas 2.2.1 and 2.2.2 respectively).
\begin{lemma}\label{hypjets}
  Let $g$ be a smooth function, the equation $X_i(f)=g$ admits a
  formal solution along the subspace $S_i$ if and only if
  $$\frac{\partial^{2k}g}{\partial x_i^k\partial y_i^k}_{\vert
    S_i}=0.$$
\end{lemma}
\begin{lemma}\label{hypflat}
  Let $g$ be a $(x_i,y_i)$-flat function along the subspace $S_i$ then
  there exists a smooth function $f$ for which $X_i(f)=g$.
\end{lemma}
{\bf Remarks:}
\begin{enumerate}
\item Let us point out that when $n=1$ the decomposition claimed in
  Proposition \ref{hyp1} had been formerly given by Guillemin and
  Schaeffer \cite{guilleminfuchsian}, by Colin de Verdi\`{e}re and Vey
  in \cite{colinvey} and Eliasson in \cite{eli1}.
\item The recipe for solving the equation specified in the lemma above
  in the case $n=1$ was given by Eliasson in \cite{eli1}. The recipe
  for the general case follows the same guidelines. It is given by the
  following formula.
\begin{equation}
f(x_1,y_1,\dots,x_n,y_n)=-\int_0^{T_i(x_1,y_1,\dots,x_n,y_n)}
g(\phi_t(x_1,y_1,\dots,x_n,y_n))dt \label{eqn:flathyp}.
\end{equation}
where $T_i$ is the function,
$$T_i(x_1,y_1,\dots,x_n,y_n)=\begin{cases}
  \frac{1}{2}\ln{\frac{x_i}{y_i}}
  \quad x_i y_i > 0 \\
  \frac{1}{2}\ln{\frac{-x_i}{y_i}} \quad x_i y_i< 0
\end{cases}$$
and $\phi_t(x_1,y_1,\dots,x_n,y_n)$ the flow of the vector field
$X_i$.  Observe that $f$ is defined outside the set
$\Omega=\Omega_1\cup \Omega_2$ where $\Omega_1$ and $\Omega_2$ are the
sets: $\Omega_1=\{(x_1,y_1,\dots,x_n,y_n), x_i=0\}$ and $\Omega_2=
\{(x_1,y_1\dots,x_n,y_n), y_i=0\}$.  In \cite{evathesis} it is proven
that $f$ admits a smooth continuation in the whole neighbourhood
considered and that it is a solution of the equation $X_i(f)=g$.
\item From the formula specified above one deduces that if $X_j(g)=0$
  for $j\neq i$ then $X_j(f)=0$.
\item In contrast to the uniqueness of the function $g_1$ in the
  decomposition obtained in proposition \ref{decomeli} for elliptic
  vector fields, the function $g_1$ specified in the decomposition is
  not unique. In fact, if $g_1$ and $h_1$ are two functions fitting in
  the decomposition their difference is an $(x_i,y_i)$ flat function
  along $S_i$. In order to check this, observe $g_1-h_1=X_i(h_2-g_2)$
  where $h_2$ is a function such that $g=h_1+X_i(h_2)$. Now, on the
  one hand the Taylor expand of $g_1-h_1$ in the $x_i,y_i$ variables
  has the form $\sum_j c_j(\check z_i)(x_i\cdot y_i)^j$ but, on the
  other hand, the Taylor expand of $X_i(h_2-g_2)$ has the form
  $\sum_{jk} c_{jk}(\check z_i)x_i^j y_i^k$ with $j\neq k$ and since
  the equality $g_1-h_1=X_i(h_2-g_2)$ holds we deduce that $g_1-h_1$
  is an $(x_i,y_i)$-flat function along $S_i$.
\item Let us show the last point of the proposition. The first step in
  the proof of the proposition was to take care of the formal Taylor
  series in $(x_i,y_i)$. Then it is easy to see that one can always
  choose Borel resummations of these formal expansions which are
  annihilated by $X_j$ ($j\neq i$) whenever $g$ is.

  Finally we integrate the flat function using the
  formula~\eqref{eqn:flathyp}, on which one can check directly that
  $f$ is invariant by the flow of $X_j$ ($j\neq i$) whenever $g$ is,
  at least in a neighbourhood of any point where the formula is well
  defined. In other words $X_j(g)=0$ implies $X_j(f)=0$ at these
  points, and hence everywhere by continuity.
\end{enumerate}

\section{A special decomposition for focus-focus vector fields}
The aim of this section is to prove the analogue to propositions
\ref{decomeli} and \ref{hyp1} for a focus-focus pair.

But before stating and proving this result we need some preliminary
material concerning the integration of equations of type $X(f)=g$ in a
neighbourhood of a hyperbolic zero (in the sense of Sternberg) of the
vector field $X$.  As we will see, the resolution of this equations is
closely related to the problem of finding the desired decomposition
for focus-focus pairs.

\subsection{Digression: Two theorems of Guillemin and Schaeffer}

A point is called a hyperbolic zero of a vector field $X$ if the
vector field vanishes at this point and all the eigenvalues of the
matrix associated to the linear part of $X$ have non-zero real part.

According to Sternberg's linearization theorem a vector field can be
linearized in a neighbourhood of a hyperbolic zero.

The following two theorems are concerned with the integration of
equations of type $X(f)=g$ in a neighbourhood of a hyperbolic zero.
These theorems A and B correspond to theorems 2 and 4 in section 4 of
\cite{guilleminfuchsian}.

\begin{thm}{\bf (Theorem A)}\cite{guilleminfuchsian}

  Let $V$ be a linear vector field on $\mathbb R^n$ with a hyperbolic
  zero at the origin and let $c$ be a fixed constant. Then given a
  smooth function $g$ flat a the origin, there exists a smooth
  function defined in a neighbourhood of the origin which is flat at
  the origin and such that:

  $$V(f)+cf=g.$$

\end{thm}

The theorem that follows is used in the proof of Theorem A. We recall
it here because we will need it in order to show the smoothness of
some constructions used in the next subsection. This theorem uses a
trick previously used by Nelson \cite{nelson} in his proof of the
Sternberg's linearization theorem.

\begin{thm}{\bf (Theorem B)}\cite{guilleminfuchsian} Let $U(t)$ be a group of linear
  transformations acting on $\mathbb R^n$. Let $N$ be a subspace of
  $\mathbb R^n$ invariant under $U(t)$ and let $E$ be the subspace of
  $\mathbb R^n$ consisting of all $x$ in $\mathbb R^n$ such that
  $$\lim_{t\rightarrow\infty}\vert\vert U(t)(x)-N\vert\vert=0.$$

  Let $g$ be a compactly supported function on $\mathbb R^n$ which is
  flat along $N$. Set

  $$f(x,s)=-\int_0^s e^{ct}g(U(t)(x))dt.$$
  Then for all multi-indices
  $\alpha$, $\lim_{s\rightarrow\infty} D^{\alpha}f(x,s)$ converges
  absolutely for all $x\in E$ and is a smooth function of $x$.
  Moreover this limit is flat along $N$.
\end{thm}

Observe that the vector field $X_i$ in a focus-focus pair $X_i,
X_{i+1}$ has a hyperbolic zero {\`a} la Sternberg on the set
$\{x_j=c_j, y_j=d_j, j\neq i, j\neq i+1 \}$ for fixed constants $c_j$
and $d_j$.

\subsection{Our proposition for focus-focus pairs}

When $i$ is the index of a focus-focus component, we denote by $S_i$
the set $S_i=\{x_i=0,y_i=0, x_{i+1}=0,y_{i+1}=0\}$. Let us state and
prove the decomposition result for focus-focus pairs.
\begin{prop}
  \label{focus-focus}
  Let $q_i,q_{i+1}$ be a focus-focus pair,
  \begin{eqnarray*}
    q_i & = & x_i y_i + x_{i+1} y_{i+1}\\
    q_{i+1} & = & x_i y_{i+1}- x_{i+1} y_i
  \end{eqnarray*}
  \noindent and let $g_1$ and $g_2$ be two functions satisfying the
  commutation relation:
  $$X_i(g_2)=X_{i+1}(g_1)$$
  \noindent Then there exists smooth  functions $f_1$, $f_2$ and $F$
  such that
  \begin{equation}
    X_j(f_k)=0 \quad j\in\{i,i+1\}\quad k\in\{1,2\}
    \label{equ:commut}
  \end{equation}
  such that
  \begin{eqnarray*}
    g_1 & = & f_1+X_i(F) \\
    g_2 & = & f_2+X_{i+1}(F)
  \end{eqnarray*}
  Moreover
  \begin{enumerate}
  \item $f_2$ is uniquely defined and satisfies $X_j(f_2)=0$ whenever
    $X_j(g_2)=0$ for some $j$;
  \item $f_1$ is uniquely modulo functions that are $z_j$-flat along
    $S_j$ and satisfy~\eqref{equ:commut};
  \item one can choose $F$ and $f_1$ such that $X_j(F)=X_j(f_1)=0$
    whenever $X_j(g_1)=X_j(g_2)=0$ for some $j\neq i$.
  \end{enumerate}
\end{prop}
{\bf Remark:} In the case $n=2$ the proposition above was proven by
Eliasson~\cite{eli1}.
\begin{proof}
  Here again the proof if a mild extension of Eliasson's.  Without
  loss of generality, one can assume that $i=1$. The flow of $X_2$
  defines an $S^1$-action which will be used in the proof. We can
  visualise this $S^1$-action easily using complex coordinates
  $z_1=x_1+ix_2$ and $z_2=y_1+iy_2$, so that $q_1+iq_2=\bar{z_1}z_2$.
  The flow of $q_2$ is the $S^1$ action given by $(z_1,z_2)\mapsto
  e^{-it}(z_1,z_2)$ whereas the flow of $q_1$ is the hyperbolic
  dynamics given by $(z_1,z_2)\mapsto (e^{-t}z_1,e^{t}z_2)$ (both
  flows act trivially on the remaining coordinates). When we say that
  a function $H$ is $S^1$-invariant for this action we mean that
  $X_2(H)=0$.

  As in the proof of Eliasson, we will first integrate along this
  $S^1$ action and then along the hyperbolic flow in an
  $S^1$-invariant way. Instead of using the formula of Eliasson (which
  consists in integrating from a transversal hyperplane through the
  origin), we will embed everything in $\RM^{2n}$ in order to apply
  the parametric versions of Theorems A and B.

  The proof consists of three steps:
  \paragraph{\textbf{1. Integrating along the $S^{1}$-action.}}

  Let $\phy_{2,t}$ be the flow of $q_2$.  As in the elliptic case
  (Proposition \ref{decomeli}) we define
  \[
  F_2=\frac{1}{2\pi}\int_0^{2\pi}(\theta-1)g_2\circ\phy_{2,\theta}
  d\theta
  \]
  and one obtains easily, by differentiating $F_2\circ\phy_{2,t}$ at
  $t=0$, that
  \begin{equation}
    X_2(F_2)=g_2-f_2,
    \label{equ:decompf2}
  \end{equation}
  where
  \begin{equation}
    f_2=\frac{1}{2\pi}\int_0^{2\pi}
    g_2(\phy_{2,\theta})d\theta,
    \label{equ:formula-f2}
  \end{equation}
  which is obviously $S^1$ invariant. Notice that if $f_2$ is any
  $S^1$ invariant function satisfying equation~\eqref{equ:decompf2}
  then by integrating along the $S^1$ flow $f_2$ is necessarily of the
  form given by~\eqref{equ:formula-f2}. Hence such an $f_2$ is indeed
  unique.

  If we check that $X_1(f_2)=0$, then we can write $g_2=f_2+X_2(F_2)$,
  with $f_2$ satisfying $X_1(f_2)=0$ and $X_2(f_2)=0$. That is to say,
  these functions $g_2$ and $f_2$ solve the second equation stated in
  the proposition.

  One can check this directly on formula~\eqref{equ:formula-f2}, using
  the commutation relation $X_1(g_2)=X_2(g_1)$ and the fact that the
  flows of $X_1$ and $X_2$ commute; one can also from
  equation~\eqref{equ:decompf2} write
  \[
  0=X_1(f_2)+X_2(X_1(F_2)-g_1),
  \]
  where $X_2(X_1(f_2))=0$. This equation can be seen as a
  decomposition for the zero function.  Using the uniqueness of the
  $S^1$-invariant function in this decomposition we obtain
  \[
  X_1(f_2)=0, \quad X_2(X_1(F_2)-g_1)=0,
  \]
  in particular this also yields that the function
  $\tilde{g_1}=g_1-X_1(F_2)$ is $S^{1}$-invariant.

  \paragraph{\textbf{2. Formal resolution of the system.}}

  In order to solve the initial system we need to find a smooth
  function $f_1$ such that $X_1(f_1)=0$ and $X_2(f_1)=0$ and a smooth
  function $F_1$ solving the system
  \begin{equation}
    \label{equ:system}
    \begin{array}{rcl}
      X_1(F_1) & = & \tilde{g_1}-f_1\\
      X_2(F_1) & = & 0,
    \end{array}
  \end{equation}
  Once this system has been solved the desired function $F$ solving
  the initial system can be written as $F=F_1+F_2$.

  In order to solve this system we will first find a formal solution
  using formal power series and in a further step we will take care of
  the remaining flat functions along $S_1$.

  We first solve the system in formal power series in $(z_1,z_2)$,
  which is fairly easy. It amounts to solving the first equation
  assuming that all terms in the series commute with $q_2$ (we can do
  this because $X_2(\tilde{g_1})=0$).  As in the hyperbolic case, the
  formal series for $f_1$ is unique and is of the form $\sum
  c_{k,\ell}(\check{z})q_1^kq_2^\ell$, where
  $\check{z}=(x_3,y_3,\dots,x_n,y_n)$.  Now we can use a Borel
  resummation in the variables $(q_1,q_2)$ for $f_1$ and an
  $S^1$-invariant Borel resummation for $F_1$, which ensures that the
  system is reduced to the situation where the right hand-side of the
  first equation of~\eqref{equ:system} is a function $g_1$ which is
  $S^1$ invariant and flat at $\{z_1=z_2=0\}$.  These Borel
  resummations can be chosen uniform in the $\check{z}$ variables.

  \paragraph{\textbf{3. Solving the equation $X_1(F_1)=g_1$ for an
      $S^1$-invariant function which is flat along $S_1$. }}

  We could finish the proof by invoking a similar formula as for the
  hyperbolic case (Proposition \ref{hypflat}). But checking the
  smoothness in all variables is not so obvious; we present here a
  small variant which uses Theorem A and B stated in the preceding
  subsection and which are contained in \cite{guilleminfuchsian}.

  The strategy is exactly the same as in \cite{guilleminfuchsian},
  with the additional requirement of keeping track of the $S^1$
  symmetry.  We give below the arguments for the sake of completeness.

  First of all, using an $S^1$-invariant cut-off function in
  $\RM^{2n}$, one can assume that $g_1$ is compactly supported while
  still commuting with $X_2$.  Again, let us call this new function by
  $g_1$. It is clear that if one solves the corresponding
  system~\eqref{equ:system} in $\RM^{2n}$, the associated germs for
  $F_1$ and $f_1$ will solve the initial local problem.  Let
  $\phy_{1,t}$ be the flow of $q_1$. The matrix associated to the
  linear vector fields $X_1$ has two positive and two negative
  eigenvalues.

  We first apply Theorem B with parameters $x_j,y_j$, $j\neq 1$ and
  $j\neq 2$ with $N=S_1$, $E=E^+=\{z_1=0\}$ and $U(t)=\phy_{1,-t}$.
  As explained in the proof of Theorem A in \cite{guilleminfuchsian}
  this allows to solve the equation to infinite order on the $2n-2$
  dimensional invariant subspace $E^+=\{z_1=0\}$.  Observe that the
  formula provided in the statement of Theorem B shows that if the
  function $g$ depends smoothly on the parameters $x_j$ and $y_j$ for
  $j\neq 1$ and $j\neq 2$ then the function $f$ does also depend
  smoothly on this parameters because $\phy_{1,-t}$ leaves the set
  $S_1$ fixed.

  Therefore using an $S^1$-invariant Borel resummation, we are then
  reduced to the case where $g_1$ is flat on $E^+$ and
  $S^1$-invariant, and we terminate by a second application of Theorem
  B with parameters $x_j,y_j$, $j\neq 1$ and $j\neq 2$ with $N=E^+$
  and $E=\mathbb R^{2n}$. That is the function $F_1$ is given by the
  formula
  \[
  F_1=-\int_0^\infty g_1\circ \phy_{1,t}dt.
  \]

  Again this function $F_1$ is smooth in all the variables since $g_1$
  is smooth in all the variables.  Using this formula we see that
  $X_2(F_1)=0$ because $\phy_{1,t}$ and $\phy_{2,\theta}$ commute.

  The justification of the last claim of the proposition goes as
  before, by examinating the explicit formulae and the Borel
  resummations. The claimed uniqueness of $f_1$ modulo $z_j$-flat
  functions along $S_j$ is a direct consequence of the uniqueness of
  the formal solution in the $z_j$ variables. Of course, one can also
  check it by an \emph{a posteriori} argument as we did in the remark
  after lemma~\ref{hypjets}.
\end{proof}
\section{The proof of Theorem~\ref{theo:principal}}
Consider $s=k_e+k_h+k_f$. As we observed in section 2. we have
$r=k_e+k_h+2k_f$. Observe also that $r=s$ if there are no focus-focus
components. We prove the theorem using induction on $s$ for a fixed
$n$.

In order to simplify the statements involving focus-focus pairs, we
introduce some more notation.  Let the vector fields
$Y_1,Y_2,\dots,Y_s$ be such that $Y_j=X_j$ for elliptic or hyperbolic
cases (\emph{ie.} for $j\leq k_e+k_h$) while
$Y_j=X_{\sigma(j)}+\sqrt{-1}X_{\sigma(j)+1}$ for focus-focus pairs
(\emph{ie.}  $j>k_e+k_h$ and $\sigma(j):=2j-k_e-k_h-1$).  Similarly we
define $\gamma_j$ to be $g_j$ for elliptic or hyperbolic indices, and
$\gamma_j=g_{\sigma(j)}+\sqrt{-1}g_{\sigma(j)+1}$ for focus-focus
indices.

For any $j\leq s$ let $\C_j$ be the space of all germs of complex
functions $f\in\Cinf(\RM^{2n},0)$ such that $Y_j(f)=\mybar{Y_j}(f)=0$,
and $\F_s=\cap_{j\leq s}\C_j$.

With these notations, the system we wish to solve has the form
$\gamma_j=f_j+Y_j(G)$ ($\forall j \in\{1,\dots,s\}$) for germs of
smooth functions $G$ and $f_j$, where $f_j\in\F_s$ and $G$ and $f_j$,
$j\leq k_h+k_e$ are real valued. The commutation relations are
$\mybar{Y_i}(\gamma_j)=Y_j(\mybar{\gamma_i})$ and
$Y_i(\gamma_j)=Y_j(\gamma_i)$ (of course the second one is redundant
except when both $Y_i$ and $Y_j$ are complex).

Suppose throughout the rest of the proof that $r<n$. For any subindex
$i$ corresponding to an elliptic or hyperbolic vector field $Y_i$ , we
denote by $z_i=(x_i,y_i)$ and
$\check{z_i}=(z_1,\dots,\check{z_i},\dots,z_{n})$. For any subindex
$j$ corresponding to a focus-focus pair $Y_j$, we denote by
$z_j=(x_i,y_i,x_{i+1},y_{i+1})$ and
$\check{z_j}=(z_1,\dots,\check{z_j},\dots,z_{n})$ (with
$i=\sigma(j)$).  We denote by $S_j$ the set $S_j=\{z_j=0\}$.

This being said, one notices that there is no more need to keep the
vector fields $Y_j$ in a particular order, which is of course most
convenient for the induction process.
\begin{sublemma}
  \label{sublemma}
  Let $Z$ be a (real or complex) vector field on $\RM^{2n}$ acting
  trivially on the variables $(z_1,\dots,z_s)$.  Let $j\leq s$. Let
  $f$ be a smooth real valued function on $\mathbb R^{2n}$ such that:
  \begin{enumerate}
  \item $f\in\F_s$
  \item $Z(f)$ is flat along $S_j$.
  \end{enumerate}
  Then there exists a smooth real valued function $\tilde{f}\in\F_r$
  such that
  \begin{enumerate}
  \item $Z(\tilde{f})=0$
  \item $f-\tilde{f}$ is flat along $S_j$.
  \end{enumerate}
\end{sublemma}
\begin{proof}
  Consider the Taylor expansion of $f$ in $z_j$. Because $Y_j(f)=0$
  this expansion is a formal series in $q_j$ (in case of an elliptic
  or hyperbolic $Y_j$) or in $q_i,q_{i+1}$ (in case of a focus-focus
  $Y_j$, with $i=\sigma(j)$). Moreover the coefficients of this
  expansion are functions of $\check{z}_j$ that are annihilated by
  $X_j$, $j\leq r, j\neq i$, and $Z$. Hence using a suitable Borel
  resummation one can come up with a smooth $\tilde{f}$ satisfying the
  requirements of our statement.
\end{proof}
\subsection{Case $s=1$}
\begin{enumerate}
\item The Cartan subalgebra has Williamson type $(1,0,0)$ or
  $(0,1,0)$.  In this case there is only one function.  Propositions
  \ref{decomeli} (in the case $X_i$ is elliptic) and \ref{hyp1} (in
  the case $X_i$ is hyperbolic) guarantee that the theorem holds.
\item The Cartan subalgebra has Williamson type $(0,0,1)$. In this
  case there are two functions $g_1$ and $g_2$ fulfilling the
  conditions specified in Proposition \ref{focus-focus}, and the
  proposition guarantees that the theorem holds.
\end{enumerate}
\subsection{ Passing from $s$ to $s+1$}
By hypothesis we can construct $G$ and $f_1,\dots,f_s$ such that
\[
\forall j\leq s, \qquad \gamma_j=f_j+Y_j(G),
\]
with $f_j\in\F_r, \forall j\leq r$.  Observe that when we pass from
$s$ to $s+1$ we are adding a real vector field if the Williamson type
changes from $(k_e,k_h,k_f)$ to $(k_e+1,k_h,k_f)$ or from
$(k_e,k_h,k_f)$ to $(k_e,k_h+1,k_f)$. In the case we increase in one
the number of focus-focus components we are adding a complex vector
field.  The proof will go in two steps. First we modify the existing
$f_j$ and $G$ in such a way that the new $f_j$'s, $j\leq s$ are in
$\F_{s+1}$.  The final step is to look for a new $G$ of the form
$\tilde{G}=G+K$ which leads to the system
\[
Y_1(K)=\dots Y_s(K)=0, \quad \tilde{\gamma}_{s+1}=f_{s+1}+Y_{s+1}(K),
\]
with $Y_j(\tilde{\gamma}_{s+1})=\mybar{Y_j}(\tilde{\gamma}_{s+1})=0$,
$\forall j\leq s$.

\paragraph{\textbf{1.}} Let us consider the commutation relations
\[
\mybar{Y_{s+1}}(\gamma_j)=Y_j(\overline{\gamma_{s+1}}) \quad \textrm{
  and } \quad Y_{s+1}(\gamma_j)=Y_j(\gamma_{s+1}).
\]
We distinguish three subcases:
\begin{enumerate}
\item The vector field $Y_j$ is elliptic: From the uniqueness of the
  function $g_1$ of the decomposition in Proposition \ref{decomeli}
  (possibly applied to the real and imaginary parts of $Y_{s+1}$) this
  condition tells us that $Y_{s+1}(f_j)=0$.  Therefore in this case no
  modification of $f_j$ is required and $f_j\in\F_{s+1}$.
\item The vector field $Y_j$ is hyperbolic: By applying lemma
  \ref{hypjets} we deduce that the $z_j$-jet of $Y_{s+1}(f_j)$ is
  zero. We can write $Y_{s+1}(f_j)=\alpha_j$ where $\alpha_j$ is a
  $z_j$-flat function along $S_j$. We can now apply sublemma
  \ref{sublemma} to obtain the following decomposition
  $f_j=\tilde{f}_j+\phi_j$ where $\tilde{f}_j\in\F_{s+1}$ and
  $\phi_j\in\F_s$ is a $z_j$-flat function.

  We may apply lemma \ref{hypflat} to the function $\phi_j$ to find a
  function $\varphi_j$ satisfying $Y_j(\varphi_j)=\phi_j$.  According
  to Proposition \ref{hyp1}, this function $\varphi_j$ can be chosen
  such that $Y_j(\varphi_j)=0$ for $j\neq i$ and $j\leq s$. Hence for
  this $\gamma_j$ we can write
  $$\gamma_j=\tilde{f}_j+Y_j(\varphi_j+G).$$
\item The vector field $Y_j$ is a focus-focus complex vector field.
  The commutation conditions also read:
  $$Y_{s+1}(\Re\gamma_j)=\Re(Y_j)(\gamma_{s+1}).$$
  $$Y_{s+1}(\Im \gamma_j)=\Im(Y_j)(\gamma_{s+1}).$$

  From the second equation and the uniqueness of the function $f_2$
  obtained in Proposition \ref{focus-focus} we obtain $Y_{s+1}(\Im
  f_j)=0$ so we only need to modify $\Re f_j$.

  Now since $\Im(Y_j)(\Re f_j)=0$ and $\Re(Y_j)(\Re f_j)=0$ we can
  invoke the uniqueness up to a flat function of the function $f_1$ in
  the decomposition claimed in proposition \ref{focus-focus} applied
  to the first equality to deduce that $Y_{s+1}(\Re f_j)$ is
  $z_j$-flat along $S_j$. Hence by sublemma \ref{sublemma} applied to
  $Z=Y_{s+1}$ we can write $\Re f_j=h_j+\phi_j$ where $h_j$ is a real
  function in $\F_{s+1}$ and $\phi_j\in\F_s$ is a real $z_j$-flat
  function along $S_j$; therefore as in the proof of
  Proposition~\ref{focus-focus} we can integrate $\phi_j$ to a
  function $\varphi_j$ satisfying $\Re Y_j(\varphi_j)=\phi_j$. Hence
  \[
  \gamma_j=\tilde{f}_j+Y_j(G+\varphi_j),
  \]
  with $\tilde{f}_j=f_j-\phi_j\in\F_{s+1}$.
\end{enumerate}
\vspace{5mm}

\paragraph{\textbf{2. }}  After considering all these cases we may
write
\[
g_j=\tilde{f}_j+Y_j(\varphi_j+G) \quad, \forall j\leq s
\]
where $\varphi_j\in\F_s$ is a real function equal to the zero function
for subindices corresponding to elliptic $Y_j$.  Now define
$\widetilde{G}=\sum_i\varphi_i+G$.  This function satisfies
\[
Y_j(\widetilde{G})= Y_j(\varphi_j+G)\quad, \forall j\leq s.
\]
To finally prove the theorem, it suffices to find a real function $K$
and $f_{s+1}\in\F_{s+1}$ such that
\[
\left\{\begin{array}{l}
    \gamma_j=\tilde{f}_j+Y_j(\tilde{G}+K), \quad \textrm{ for } j\leq s\\
    \gamma_{s+1}=f_{s+1}+Y_{s+1}(\tilde{G}+K).
\end{array}\right.
\]
But consider $\tilde{\gamma}_{s+1}:=\gamma_{s+1}-Y_{s+1}(\widetilde
G)$.  The commutation relations yield
\begin{equation}
  Y_j(\tilde{\gamma}_{s+1})=\mybar{Y_j}(\tilde{\gamma}_{s+1})=0
  \label{equ:gamm-comm}
\end{equation}
for $j\leq s$, and we still have (in case $s+1$ is a focus-focus
index)
\begin{equation}
  Y_{s+1}(\mybar{\tilde{\gamma}_{s+1}})=\mybar{Y_{s+1}}(\tilde{\gamma}_{s+1}).
  \label{equ:ff-comm}
\end{equation}
Thus our system becomes
\[
\left\{\begin{array}{l}
    0=Y_j(K), \quad \textrm{ for } j\leq s\\
    \tilde{\gamma}_{s+1}=f_{s+1}+Y_{s+1}(K),
\end{array}\right.
\]
and since $\tilde{\gamma}_{s+1}\in\F_s$
(equation~\eqref{equ:gamm-comm}), it is solved by an application of
proposition~\ref{decomeli}, \ref{hyp1} or \ref{focus-focus}, depending
on the type of $Y_{s+1}$ (notice that the relation~\eqref{equ:ff-comm}
is precisely the commutation relation required in the focus-focus
case).  This ends the proof of the theorem.

\section{Deformations of completely integrable systems}
\label{sec:deformation}

Theorem~\ref{theo:principal} has a natural interpretation in terms
of infinitesimal deformations of integrable systems near
non-degenerate singularities.  This was stated without proof in
\cite{san-habil}. Let us recall briefly the appropriate setting.

A completely integrable system on a symplectic manifold $M$ of
dimension $2n$ is the data of $n$ functions $f_1,\dots,f_n$ which
pairwise commute for the symplectic Poisson bracket: $\{f_i,f_j\}=0$
and whose differentials are almost everywhere linearly independent.

When we are interested in geometric properties of such systems,
the main object under consideration is the (singular) lagrangian
foliation given by the level sets of the momentum map
$f=(f_1,\dots,f_n)$.  We introduce the notation $\mathbf{f}$ for
the linear span (over $\RM$) of $f_1,\dots,f_n$. It is an
$n$-dimensional vector space. It is also an abelian Poisson
subalgebra of the Poisson algebra $X=(\Cinf,\{,\})$. Let
$\C_\mathbf{f}=\{h\in X, \{\mathbf{f},h\}=0\}$ be the set of
functions that commute with all $f_i$.  By Jacobi identity
$\C_\mathbf{f}$ is a Lie subalgebra of $X$. The fact that
$df_1\wedge\cdots\wedge df_n\neq 0$ almost everywhere implies that
$\C_\mathbf{f}$ is actually abelian. From now on, we are given a
point $m\in M$ and everything is localised at $m$; in particular
$X$ is the algebra of germs of smooth functions at $m$.

\begin{defn}
  Two completely integrable systems $\mathbf{f}=\gener{f_1,\dots,f_n}$
  and $\mathbf{g}=\gener{g_1,\dots,g_n}$ are \emph{equivalent} (near
  $m$) if and only if
\[
\C_\mathbf{f} = \C_\mathbf{g}
\]
\end{defn}
Geometrically speaking, $\mathbf{f}$ is equivalent to $\mathbf{g}$ if
and only if the functions $f_i$ are constant along the leaves of the
$\mathbf{g}$-foliation (or vice-versa).

We wish to describe infinitesimal deformations of integrable systems
modulo this equivalence relation. For this we fix an integrable system
$\mathbf{f}$ and introduce a deformation complex as follows.  Let
$L_0\simeq \RM^n$ be the typical commutative Lie algebra of dimension
$n$.  $L_0$ acts on $X$ by the adjoint representation:
\[
L_0\times X\ni(\ell,g) \mapsto \{\mathbf{f}(\ell),g\}\in X .
\]
Hence $X$ is an $L_0$-module, in the Lie algebra sense, and we can
introduce the corresponding Chevalley-Eilenberg
complex~\cite{chevalley-eilenberg}: for $q\in\NM$,
$C^q(L_0,X)=\textup{Hom}({L_0}^{\wedge q},X)$ is the space of
alternating $q$-linear maps from $L_0$ to $X$ (regarded merely as
real vector spaces), with the convention $C^0(L_0,X)=X$. The
associated differential is denoted by $d_\mathbf{f}$. Following
\cite{chevalley-eilenberg} for a $0$-cochain  $g\in X$, the
$1$-cochain $d_\mathbf{f}(g)$ is  $d_f (g)(l)
=\{\mathbf{f}(l),g\}, l \in L$ and for a $k$-cochain  $\phi$ the
$k+1$ cochain $d_\mathbf{f}(\phi)$ is
$$d_\mathbf{f}(\phi)(l_1,\dots,l_{k+1})=
\frac{1}{k+1}\sum_{i=1}^{k+1}(-1)^{i+1}\{\mathbf{f}(l_i),\phi(\check{l_i})\},\l_i
\in L,$$

\noindent where $\check{l_i}=(l_1,\dots,\check{l_i},\dots, l_{k+1}).$

Now since $L_0$ acts trivially on $\O_\mathbf{f}$, the quotient Lie
algebra $X/\O_\mathbf{f}$ is a $L_0$-module, and we can define the
corresponding Chevalley-Eilenberg complex: for $q\in\NM$,
$C^q(L_0,X/\O_\mathbf{f})=\textup{Hom}({L_0}^{\wedge
  q},X/\O_\mathbf{f})$, with differential denoted by
$\bar{d}_\mathbf{f}$.

Finally we define the \emph{deformation complex}
$C^\bullet(\mathbf{f})$ as follows:
\[
\xymatrix{0 \ar[r] &%
 X/\O_{\mathbf{f}}  \ar[r]^-{\bar{d}_{\mathbf{f}}} &%
 C^1(L_0,X/\O_{\mathbf{f}}) \ar[r]^-{\partial_{\mathbf{f}}} &%
C^2(L_0,X) \ar[r]^-{d_{\mathbf{f}}} &%
C^3(L_0,X) \ar[r]^-{d_{\mathbf{f}}} &%
\cdots}
\]
where $\partial_{\mathbf{f}}$ is defined by the following diagram,
where all small triangles are commutative
($C^k(L_0,\O_{\mathbf{f}})$ is always in the kernel of
$d_{\mathbf{f}}$)~:
\[
\xymatrix{ 0 \ar[r] &%
  X \ar[r]^-{d_{\mathbf{f}}} \ar[d]_{\pi} &%
  C^1(L_0,X) \ar[r]^-{d_{\mathbf{f}}} \ar[d]_{\pi} &%
  C^2(L_0,X) \ar[r]^-{d_{\mathbf{f}}} \ar[d]_{\pi} &%
  {\dots} \\
  0 \ar[r] &%
  X/\O_{\mathbf{f}} \ar[r]_-{\bar{d}_{\mathbf{f}}}
  \ar[ur]^{\partial_{\mathbf{f}}} &%
  C^1(L_0,X/\O_{\mathbf{f}}) \ar[r]_-{\bar{d}_{\mathbf{f}}}
  \ar[ur]^{\partial_{\mathbf{f}}}&%
  C^2(L_0,X/\O_{\mathbf{f}}) \ar[r]_-{\bar{d}_{\mathbf{f}}}
  \ar[ur]^{\partial_{\mathbf{f}}} &%
  {\dots} }
\]
For all cochain complexes, cocycles and coboundaries are denoted
the standard way: $Z^q(\cdot)$ and $B^q(\cdot)$. In the analytic
category a similar deformation complex was introduced recently by
Van Straten and Garay (\cite{garaystraten} and \cite{garay}) and
(for the first degrees) by Stolovitch \cite{stolovitch-ihes}. The
equivalence used in the analytic category is much easier to handle
due to the absence of flat functions.

\begin{defn}
  $Z^1(\mathbf{f})$ is the space of infinitesimal deformations of
  $\mathbf{f}$ modulo equivalence.
\end{defn}
If we fix a basis $(e_1,\dots,e_n)$ of $L_0$, a cocycle $\alpha\in
Z^1(\mathbf{f})$ is just a set of functions
$g_1=\alpha(e_1),\dots,g_n=\alpha(e_n)$ (defined modulo
$\O_\mathbf{f}$) such that
\begin{equation}
  \label{equ:cocycle}
  \forall i,j\qquad \{g_i,f_j\} = \{g_j,f_i\} .
\end{equation}
It is an infinitesimal deformation of $\mathbf{f}$ in the sense that,
modulo $\epsilon^2$,
\[
\{f_i+\epsilon g_i,f_j+\epsilon g_j\} \equiv 0 .
\]

A special type of infinitesimal deformations of $\mathbf{f}$ is
obtained by the infinitesimal action of the group $G$ of local
symplectomorphisms: given a  function $h\in X$ one can define the
deformation cocycle $\alpha\in Z^1(\mathbf{f})$ by
\begin{equation}
  L_0\ni\ell\mapsto
  \alpha(\ell)=\{h,\mathbf{f}(\ell)\} \mod \O_\mathbf{f}.
  \label{equ:coboundary}
\end{equation}
In other words, the set of all such cocyles, with $h$ varying in
$X$, is the orbit of $\mathbf{f}$ under the adjoint action on
$Z^1(\mathbf{f})$ of the Lie algebra of $G$. From
equation~\eqref{equ:coboundary} one immediately sees that this
orbit is exactly $B^1(\mathbf{f})$.

In the particular case that $\omega$ is the Darboux symplectic
form $\omega_0=\sum_{i=1}^n dx_i\wedge dy_i$ and
$\mathbf{f}=(q_1,\dots,q_n)$ is a Williamson basis as specified in
theorem \ref{Willi}, we can reformulate the statement of
theorem~\ref{theo:principal} in cohomological terms.

 Namely, in
this case since $\{f_i,f\}= X_i(f)$,  we can write $C_q=\{ f\in X, X_i(f)=0,
\forall i\}$. Let $\alpha$ be a 1-cocycle, the cocycle condition specified in
formula \ref{equ:cocycle}  reads as $X_j(g_i)=X_i(g_j)$ where
$g_i=\alpha(e_i)$. But this is nothing but the commutation hypothesis of
theorem~\ref{theo:principal} therefore  there exists a function $G$ such that
$g_i=f_i+ X_i(G)$. Using formula \ref{equ:coboundary} and the definition of
$g_i$ this shows that $\alpha$ is a coboundary. In other words, what theorem
2.2 shows in cohomological terms is that any $\alpha\in Z^1(\mathbf{f})$ is
indeed a coboundary.  And this proves the following reformulation of
theorem~\ref{theo:principal}:
\begin{thm}
  Let $q_1,\dots,q_n$ be a standard basis (in the sense of Williamson)
  of a Cartan subalgebra of $\mathcal{Q}(2n,\RM)$. Then the
  corresponding completely integrable system $\mathbf{q}$ in
  $\RM^{2n}$ is $\Cinf$-infinitesimally stable at $m=0$:
  that is,
  \[
  H^1(\mathbf{q})=0 .
  \]
\end{thm}

\begin{rem}
  Our proof actually shows that the result is also true when we
  include a smooth dependence on parameters in the definition of the
  deformation complex.
\end{rem}

This theorem should have important applications in semi-classical
analysis, where we consider pseudodifferential operators with
$\Cinf$ symbols depending on a small parameter $\hbar$.
One can define a similar deformation complex for pseudodifferential
operators, where the deformation is understood with respect to the
parameter $\hbar$.  Then in many situations the vanishing of the
classical $H^1$ implies the vanishing of the pseudodifferential $H^1$.
See \cite{san-habil} for general remarks, and
\cite{san-focus,colin-singularites} for applications in simple cases
where the vanishing of the pseudodifferential $H^1$ was checked
explicitly.

\bibliographystyle{amsplain} \providecommand{\bysame}{\leavevmode\hbox
  to3em{\hrulefill}\thinspace}

\end{document}